\documentclass[12pt,a4paper]{extarticle}
\usepackage[english]{babel}
\usepackage{amsmath, amsthm, amsfonts, amssymb, bbm}
\usepackage{mathtools}

\usepackage{geometry}
\newcommand{\ov}{\overline}

\newcommand{\vf}{\varphi}
\newcommand{\ve}{\varepsilon}
\newcommand{\cF}{{\mathcal F}}

\newcommand{\cL}{{\mathcal L}}
\newcommand{\cD}{{\mathcal D}}
\newcommand{\1}{\pmb{1}}
\newcommand{\mbR}{{\mathbb R}}
\newcommand{\mbN}{{\mathbb N}}

\newcommand{\sign}{\mathop{\rm sign}}
\newcommand{\e}{\mathrm{e}}

\newtheorem{theorem}{Theorem}[section]

\newtheorem{proposition}[theorem]{Proposition}
\newtheorem{corollary}[theorem]{Corollary}
\theoremstyle{definition}
\newtheorem{definition}[theorem]{Definition}

\theoremstyle{remark}
\newtheorem{remark}[theorem]{Remark}
\numberwithin{equation}{section}
\begin{document}
\title{Arratia flow with drift and Trotter formula for Brownian web}

\author{A.A.Dorogovtsev, M.B.Vovchanskii}



\maketitle

\allowdisplaybreaks

\begin{abstract}
An analog of the Trotter formula for the Arratia flow is presented. Perturbations of the Brownian web by mappings associated with an ordinary differential equation with a smooth right part are considered and proved to be convergent exclusively in the weak sense. The flow obtained as a limit is the Arratia flow with drift. 
\end{abstract}

\section{Introduction}
\label{intro} 
In this article the fractional step method for the Brownian web is developed. We consider perturbations of the web by an external force. The resulting motion is obtained as a limit of approximations where the web and a smooth flow of diffeomorphisms are applied subsequently on small time intervals. Compared to the case when the flows are generated by stochastic differential equations with smooth coefficients we prove only the weak convergence instead of the convergence in mean or in probability. It is an essential feature of the Brownian web. We show that the sum of its increments in some point converges only in the weak sense.

We start with the definition of the Brownian web. Let $(\cF_t)_{t\in\mbR}$ be a filtration on some probability space, and  $(u_1, t_1), \ldots, (u_N, t_N) \in\mbR^2.$ For an $N-$tuple of continuous random processes  $\{B_{t_1, t}(u_1)\mid t\geq t_1\}, \ldots,$ \\ $\{B_{t_N, t}(u_N)\mid t\geq t_N\}$ define 
$$
\theta_{ij}=\inf\left\{s\geq t_i\vee t_j \mid B_{t_i, s}(u_i)=B_{t_j, s}(u_j)\right\}.
$$
Suppose $B_{t_k, t_k}(u_k)=u_k, k=\overline{1, N}.$
\begin{definition}
\label{coalescing.bm}
{\sloppypar
The $N-$tuple of continuous random processes $\{B_{t_1, t}(u_1)\mid t \geq t_1\}, \ldots,$ $\{B_{t_N, t}(u_N)\mid t\geq t_N\}$ is called coalescing Brownian motions starting from $u_k$ at time $t_k, k =\overline{1, N},$ w.r.t. the filtration $(\cF_t)_{t\in\mbR}$ if $\{B_{t_j, t}(u_j)\mid t\geq t_j\}$ is a $(\cF_t)_{t\ge t_j}$-martingale, and
$\{B_{t_i, t}(u_i)B_{t_j, t}(u_j)-(t-\theta_{ij})_{+}\mid t\geq t_i\vee t_j$\} is a $(\cF_t)_{t\geq t_i\vee t_j}$-martingale, $i, j=\overline{1,N}.$ \par}
\end{definition}

\begin{definition}
\label{defn1}
{
\sloppypar
A Brownian web (see, for instance, \cite{BW:convergence,BW:full,dynamics}) is a collection of random processes $\left\{\vf_{t,\cdot}(u) \in C([t;+\infty))\mid u,t\in\mbR \right\}$ such that, 
given $(u_1, t_1), \ldots, (u_N, t_N)$ the processes
$\vf_{t_1, \cdot}(u_1), \ldots, \vf_{t_N, \cdot}(u_N)$
are coalescing Brownian motions in the sense of Definition~\ref{coalescing.bm} w.r.t. the filtration
$$
\cF_t=\sigma\left(\vf_{s, r}(u), s\leq r\leq t, u\in\mbR\right), t\in\mbR.
$$
\par}
\end{definition}

The filtration $(\cF_t)_{t\in\mbR}$ may be referred to as a one generated by the Brownian web, and $\vf_{s, t}(u)$ may be interpreted as a position at time $t$ of a particle that starts from $u$ at time $s$ and is carried by the web.  
\begin{remark}
The family $\Phi$ of random variables $\{\vf_{s, t}(u)\mid s\leq t\,, u\in\mbR\}$ can be treated as a family of random mappings: for any pair $(s, t)$ $u\mapsto \vf_{s, t}(u)$ is a mapping from $\mbR$ into $\mbR.$ Due to the coalescence property of the Brownian web these mappings are discontinuous yet monotone increasing. It can be proved~\cite{dorogovtsev:continuity} that there exists a modification of $\Phi$ such that $\vf_{s, t}\colon\mbR\mapsto\mbR$ is a rcll function. As a result, one can consider $\Phi$ as a family of measurable rcll mappings from $\mbR$ into $\mbR$ indexed by pairs of time marks. The following property, similar to those of stochastic flows, is observed: for any fixed $p\leq q\leq r$ and any $u,$ with probability $1$~\cite{dynamics}
\begin{equation}
\label{eq1}
\vf_{q, r}(\vf_{p, q}(u))=\vf_{p, r}(u).
\end{equation}

The Brownian web can be considered as an example of a stochastic dynamical system in the sence of~\cite{arnold}. It is shown in~\cite{riabov} that one can construct a modification of the Brownian web in such a way that it becomes a right cocycle.
However, the proposed modification does not possess the rcll property stated above.  Another facts about properties of different modifications of the Brownian web are discussed in~\cite{BW:convergence,BW:full}.
\end{remark}

Consider a family of random processes $\{Y(u)=\vf_{0,\cdot}(u) \in C([0;\infty])\mid u\in\mbR\}.$ This object, called the Arratia flow, can be understood  as a set of Brownian particles that start from all points of the real axe simultaneously at time $0$ and are carried by the web. Each two of them move independently before a collision and merge after it. Note that the Arratia flow was introduced in~\cite{arratia, selfrepelling} as a limit of rescaled random walks. From the properties of the Brownian web it follows   
\begin{enumerate}
\item for any $u$ $Y_\cdot(u)$ is a Brownian motion w.r.t. the joint filtration;
\item for $u_1\leq u_2$ \ $Y(u_1)\leq Y(u_2);$
\item for any $u_1, u_2$\quad 
\begin{equation*}
\label{arratia.flow.char}
\left\langle Y(u_1), Y(u_2) \right\rangle_t =\left( t-\inf{\{r\mid Y_r(u_1)=Y_r(u_2)\}} \right)_{+} = \int_0^t \1_{\{ Y_r(u_1)=Y_r(u_2)\}} dr .
\end{equation*}
\end{enumerate}

In order to allow a more general law of the motion of particles inside the Arratia flow one can consider Brownian motions with drift. Such non-zero term of finite variation introduces an external mechanical force into the picture. This leads to the following definition (taken from~\cite{dorogovtsev:mono} after some reformulation). 
\begin{definition}
\label{defn2}
Let $a$ be a measurable function on $\mbR.$ A family of random processes $\{Y^a(u)|u\in\mbR\}$ is called the Arratia flow with a drift $a$ if
\begin{enumerate}
 \item  for any $u$ 
$$
Y^a_t(u)=u+\int_0^t a(Y^a_s(u))ds+ B_t(u),
$$
where $B(u)$ is a Brownian motion w.r.t. the filtration $\cF^{Y^a}$ generated by $\{Y^a(u)|u\in\mbR\},$
\item for any $u_1, u_2$\quad 
\begin{equation*}
\left\langle B(u_1), B(u_2) \right\rangle_t =\left( t-\inf{\{r\mid Y^a_r(u_1)=Y^a_r(u_2)\}} \right)_{+} = \int_0^t \1_{\{ Y^a_r(u_1)=Y^a_r(u_2)\}} dr .
\end{equation*}
\end{enumerate}
\end{definition}

The existence of the Arratia flow with a bounded $a$ satisfying the Lipschitz condition is proved in~\cite{dorogovtsev:mono}. The same proof admits an extension to the case of an unbounded Lipschitz continuous drift.

Definition~\ref{defn2} guarantee that any two processes $Y^a(u_1)$ and $Y^a(u_2)$ coalesce after the meeting and, roughly speaking, are independent before the meeting occurs. The reasons of the coalescence appearing can be briefly explained as follows. Firstly, note  that the moment $\theta=\inf{\{r\mid Y^a_r(u_1)=Y^a_r(u_2)\}}$ is a Markov moment w.r.t. $\cF^{Y^a}.$ Secondly, for $k= 1,2$,  the process $Y^a(u_k)$ solves a stochastic differential equation $dY^a_t(u_k) = a(Y^a_t(u_k))dt + dB_t(u_k)$ with some Brownian motion $B(u_k), k = 1,2.$ Due to Property 2, the Brownian motions $B(u_1)$ and $B(u_2)$ coincide after the moment $\theta.$  Either process $t\mapsto Y^a_{t+\theta}(u_k), k = 1,2,$ still solves the stochastic differential equation $dz_t = a(z_t)dt + dB_{t+\theta}(u_1)$ with $z_0=Y^a_\theta(u_1)=Y^2_\theta(u_2).$ For a ''good enough'' drift coefficient such equation has a unique strong solution. Thus $Y^a(u_1)=Y^a(u_2)$ after the moment~$\theta.$

The main goal of the paper is to obtain an Arratia flow with drift from a Brownian web via perturbations of a Brownian web by the flow of solutions to a deterministic equation $dz_t = a(z_t)dt.$ More precisely, one can expect that the action of the Brownian web and this flow on small time intervals subsequently and repeatedly allows to incorporate the drift into the Arratia flow in a way compared to that of the well-known Trotter formula~\cite{reed}. The Trotter formula, postulated for flows driven by vector fields, either deterministic or stochastic ones, states that it is possible to decompose an external influence of a finite number of forces into a sum of their actions on subsequent intervals of time. In the case of stochastic flows driven by a SDE with smooth coefficients results of such kind can be found in \cite{kotelenez}. We propose a method to build perturbations of a Brownian web by a deterministic vector field $a.$ This means that the flow of solutions to $dz_t = a(z_t)dt$ and the mappings  $\{\vf_{s, t}\mid s\leq t\}$ are applied in turns. For a detailed and precise formulation, refer to Equations~(\ref{eq1'})-(\ref{eq2'}) and definitions there. The weak convergence of $N$-point motions of the perturbed web to those of an Arratia flow with a drift $a$ is established (Theorem \ref{thm1}).

Before proceeding further we want to discuss the property of the Brownian web that illustrates the complexity of a noise associated with it and, at the same time, points out the difference compared to the deterministic Trotter formula, namely, that the web is not a flow of solutions to some ``good'' SDE. For a usual SDE
\begin{equation}
\label{eq2}
dX_t=\sigma(X_t)dw_t,
\end{equation}
with $w$ being a standard Brownian motion started at 0, $\sigma\in C^2(\mbR),$ the expression $\sigma(u) dw_t$ may be interpreted as an ``infinitesimal vector field'' that defines increments of the trajectory of $X$. This field contains information about the Brownian motion $w$ and the diffusion coefficient $\sigma$ and can be restored from observations of the flow. To see that, consider $X_{s, \cdot}(u),$  a solution to \eqref{eq2} started from $u$ at time $s.$ Then $\{X_{s, t}(\cdot)\mid s\leq t\}$ is a flow of diffeomorphisms~\cite{kunita}, and one may check that, for any $t,$
$$
\sum^{n-1}_{k=0}\Big(X_{t\frac{k}{n}, t\frac{k+1}{n}}(u)-u\Big)\underset{n\to\infty}{\overset{P}{\longrightarrow}}\int^t_0\sigma(u)dw_s=\sigma(u)w_t.
$$
The mappings $\{\vf_{s, t}| s\leq t\}$  associated with a Brownian web $\vf$ demonstrate a significantly different behaviour. The next proposition serves as an illustrative example.
\begin{proposition}
 \label{bw:nonconvergence}
 For any $u$ the sequence $\left\{\sum\limits_{k=0}^{n-1}(\vf_{\frac{k}{n},\frac{k+1}{n}}(u)-u)\mid n\ge 1\right\}$ does not converge in probability. 
\end{proposition}

Proof. We need use analogs of some results obtained in~\cite{dorogovtsev:density} for the Arratia flow, extending them to the case of the Brownian web. To avoid a departure from the central topic of the paper we omit minor technical details. 

Given a Brownian web $\{\vf_{t,s}(u)\mid t\le t, u\in\mbR \}$ consider an embedded Arratia flow $\left\{Y(u) = \vf_{0,\cdot}(u)\mid u\in\mbR\right\}.$ Fix $U > 0.$ Let $\mathcal{U}=(u_1, \cdot, u_N)$ be an ordered subset of $[0;U].$  For a function $\alpha=(\alpha_1,\ldots,\alpha_N)\in (C([0;1]))^N$ define
\begin{align*}
\mathcal{I}^N_{s,t}(\mathcal{U};\alpha) & = \sum_{k=1}^N\int\limits_s^{\tau_k(\mathcal{U})}\alpha_k(r)d\vf_{s, r}(u_k), \\
\mathcal{J}^N_{s,t}(\mathcal{U};\alpha) & = \sum_{k=1}^N\int\limits_s^{\tau_k(\mathcal{U})}\alpha_k^2(r)dr,
\end{align*}
where 
$$
\begin{cases}
\tau_1(\mathcal{U}) = t, \\
\tau_k(\mathcal{U}) = \inf\left\{t; r\mid \prod_{j = \overline{1, k-1}}(\vf_{s, r}(u_{j}) - \vf_{s, r}(u_{k}))\right\}, k =\overline{2, N}.
\end{cases}
$$
Define 
$$
\mathcal{E}^N_{s, t}(\mathcal{U};\alpha)=\exp(\mathcal{I}^N_{s,t}(\mathcal{U};\alpha)-\frac{1}{2}\mathcal{J}^N_{s,t}(\mathcal{U};\alpha)).
$$ 

It is shown in~\cite{dorogovtsev:density} that, for any $\{u_n\}_{n\ge 1}$ dense in $[0;U],$ a set
$$
\left\{  \mathcal{E}^N_{0, t}(\{u_1,\ldots,u_N\};\alpha)\mid \alpha\in (C([0;1]))^N, N\in\mbN\right\} 
$$ is a total set in $L^2\left(\sigma(Y_{s, r}(u), 0\leq s\leq r\leq t, u\in [0;U]\right).$ Given
$$
\cF^{U}_{s, t}=\sigma\left(\vf_{r_1, r_2}(u), s\leq r_1\leq r_2\leq t, u\in [0;U]\right), s\leq t, 
$$
 the same technique applied in~\cite{dorogovtsev:density} allows to extend this conclusion to the case of $L^2\left(\cF^{U}_{0,t}\right),$ though  we need a weaker statement that can be formulated as follows. Suppose $\xi\in L^2\left(\cF^{U}_{0,t}\right).$ If for any $N\in\mbN,\alpha\in (C([0;1]))^N, s, t\in[0;1]$ and any ordered subset $\{u_1,\ldots, u_N\}$ of $[0;U]$
$$
E\xi  \mathcal{E}^N_{0, t}(\{u_1, \ldots, u_N\};\alpha) = 0,
$$ then $\xi$ is $0$ a.s..

Thus the statement of Proposition~\ref{bw:nonconvergence} is a consequence of:
\begin{proposition}
\label{tech}
Let $\vf$ be a Brownian web. Then for any $0\leq s\leq t\leq 1, \alpha\in (C([0;1]))^N$ and any ordered set $\mathcal{U}$ of size $N$
$$
E\sum^{n-1}_{k=0}\vf_{\frac{k}{n}, \frac{k+1}{n}}(0)\mathcal{E}^N_{s, t}(\mathcal{U};\alpha)\to0, \ n\to\infty.
$$
\end{proposition}
Indeed, the variance of $\sum^{n-1}_{k=0}\vf_{\frac{k}{n}, \frac{k+1}{n}}(0)$ can be easily seen to be $1,$ so we have a contradiction.

Being rather technical, the proof of Proposition~\ref{tech} is put in Appendix. 
\qed
 
The statement of Proposition~\ref{bw:nonconvergence} can be interpreted as that the ``infinitesimal random field'' that drives the motion of particles inside the web does not allow restoration from observations of the web. This result relates to the theory of noises associated with families of mappings \cite{tsirelson,tanaka, watanabe.warren, web.black.2d}, which states the noise of the Brownian web to be ``black''.  

As it has been mentioned before, the family of mappings $\{\vf_{s,t}\mid s\leq t\}$ of the Brownian web can be treated as a dynamical system on $\mbR.$ The previous statement means that this family is not generated by a ``good'' infinitesimal vector field. But it occurs that we can still consider perturbations of this non-existing field by using the fractional step method. 

\section{Fractional Step Method for Brownian Web}
\label{Trotter_formula}
Through the whole paper, a function $a$ satisfies the Lipschitz condition on the real axe with a constant $C_a,$ and $\{A_t(u) | t\geq0\},$ for any $u,$ is a solution to the Cauchy problem
$$
\begin{cases}
dA_t(u)=a(A_t(u))dt, \\
A_0(u)=u.
\end{cases}
$$

Consider a sequence of partitions of $[0;1], \  \{t^{(n)}_0, \ldots, t^{(n)}_{N^{(n)}}\}, n\in\mbN,$ where 
$$
\lambda^{(n)}=\max_{k=\ov{0, N^{(n)}-1}}(t^{(n)}_{k+1}-t^{(n)}_{k})
$$
tends to zero as $n$ grows to infinity. For fixed $n$ and $k\in\ov{0, N^{(n)}-1}$ and $t\in[t^{(n)}_k; t^{(n)}_{k+1})$ define
\begin{align}
X^{(n)}_t(u) & =\vf_{t^{(n)}_k, t}\left(\mathop{\circ}^k_{j=1}A_{t^{(n)}_j-t^{(n)}_{j-1}}(\vf_{t^{(n)}_{j-1}, t^{(n)}_j}(\cdot))\right)(u), \label{eq1'} \\
\Delta^{(n)}_k(u) & =X^{(n)}_{t^{(n)}_k}(u)-X^{(n)}_{t^{(n)}_k-}(u), \label{eq2'}
\end{align}
where $X^{(n)}_{t^{(n)}_k-}(u)=\lim_{s\nearrow t^{(n)}_k}X^{(n)}_s(u)$ and the $\sign\circ$ stands for a composition of functions. Put $X^{(n)}_1(u)$ to equal $X^{(n)}_{1-}(u).$ 

From now we suppose that all considered $\sigma$-fields have been completed and augmented in a usual way. 

\begin{proposition}
\label{prop2'}
For all $n$ and $k$ and every $u$ a.s.
$$
|\Delta^{(n)}_k(u)|\leq\sup_{t\in[0; 1]}|a(X^{(n)}_t(u))|\cdot e^{C_a\lambda^{(n)}}(t^{(n)}_{k}-t^{(n)}_{k-1}).
$$
\end{proposition}
Proof.
From the definition of $X^{(n)}$ it follows that
$$
X^{(n)}_{t^{(n)}_k}(u)=A_{t^{(n)}_{k}-t^{(n)}_{k-1}}(X^{(n)}_{t^{(n)}_k-}(u)).
$$
An application of the Gronwall--Belmann lemma easily implies what is stated.
\qed

Define
\begin{align*}
A^{(n)}_t(u) & = \sum_{k: t^{(n)}_k\leq t}\Delta^{(n)}_k(u), \ t\in[0; 1], n\in \mbN, \\
m^{(n)}_t(u) & =X^{(n)}_t(u)-A^{(n)}_t(u), \ t\in[0; 1], n\in \mbN. 
\end{align*}

\begin{proposition}
\label{prop:characteristic_of_diff_m}
For any $u$ and $n$ \ $m^{(n)}(u)$ is a Brownian motion started at $u$ with respect to the filtration $(\cF_t)_{t\in[0; 1]},$ and, for any $u_1, u_2,$
$$
\langle m^{(n)}(u_1), m^{(n)}(u_2)\rangle_{t} = \int_0^t \1_{\{X^{(n)}_s (u_1) = X^{(n)}_s(u_2)\}} ds = \big( t - \tau^{(n)}\big)_{+},
$$
where $\tau^{(n)} = \inf\{1; s\mid X^{(n)}_s(u_1) =  X^{(n)}_s(u_2)\}.$
\end{proposition} 
Proof. One can check in a standard way that $m^{(n)}(u)$ is a continuous $(\cF_t)_{t\in[0; 1]}$-martingale with the characteristic $t$ and then apply the Levy theorem on characterization of the Brownian motion, as the $m^{(n)}$ have been built continuous.

Note that the definition of the Brownian web implies the following. For any $t^{(n)}_{j-1}\le s\le t<t^{(n)}_{j}, j = \overline{1,N^{(n)}}$ and any $v$ 
$$
X^{(n)}_t(v) = \vf_{s, t}(X^{(n)}_s(v)).
$$
Also for same $s,t$ and any $u_1, u_2,$
\begin{align*}
& E\left( X^{(n)}_t(u_1) X^{(n)}_t(u_2) - \!\int_s^t \1_{\{X^{(n)}_r (u_1) = X^{(n)}_r(u_2)\}} dr \left\vert \cF^{BW}_s \right.\right) \\
& = 
E\left( \vf_{s, t}(X^{(n)}_s(u_1)) \vf_{s, t}(X^{(n)}_s(u_2)) - \!\int_s^t \!\!\1_{\{\vf_{s, r}(X^{(n)}_s(u_1))=\vf_{s, r}(X^{(n)}_s(u_2))\}} dr \left\vert \cF^{BW}_s \right. \right) \\
& = E\left( \vf_{s, t}(v_1) \vf_{s, t}(v_2)\! - \!\int_s^t \!\!\1_{\{\vf_{s, r}(v_1))=\vf_{s, r}(v_2))\}} dr \left\vert \cF^{BW}_s\right. \right)\!\Bigg|_{v_1 = X^{(n)}_s(u_1), v_2 = X^{(n)}_s(u_2)}  
\\
& = v_1 v_2\Big|_{v_1 = X^{(n)}_s(u_1), v_2 = X^{(n)}_s(u_2)} = X^{(n)}_s(u_1) X^{(n)}_s(u_2),
\end{align*}
since, for $t\in[s; t^{(n)}_j), v_1, v_2 \in\mbR,$
$$
\langle \vf_{s, \cdot}(v_1)), \vf_{s, \cdot}(v_2) \rangle_{t} = \int_s^t \1_{\{\vf_{s, r}(v_1))=\vf_{s, r}(v_2))\}} dr.
$$

If $s,t$ belong to different intervals of the partition $\{t^{(n)}_0, \ldots, t^{(n)}_{N^{(n)}}\}$ one should carefully apply a standard reasoning based on  conditioning and two equalities just mentioned. We omit the details.

{\sloppypar
To prove the second equality of the proposition note that the difference $X^{(n)}(u_2)-X^{(n)}(u_1)$ never leaves 0 after hitting it, the latter being a consequence of the Brownian web's definition. Thus 
$$
\int_0^t \1_{\{X^{(n)}_s (u_1) = X^{(n)}_s(u_2)\}} ds = \big( t - \tau^{(n)}\big)_{+}.
$$
\par}
\qed

\section{Existence and Characterization of Weak Limits of $X^{(n)}$}
The $N$-point motions of the perturbed Brownian web $(X^{(n)}(u_1), \ldots,$ $X^{(n)}(u_N))$ are considered in $\cD([0; 1]),$ the Skorokhod space of rcll functions on $[0; 1]$ endowed with the distance~\cite{billingsley} 
\begin{align*}
d(f, g)= \inf \{\ve>0|\exists\lambda\in\Lambda: \sup_{t\in[0; 1]}|f(t)-g(\lambda(t))|  \vee\sup_{\begin{subarray}{c}
t,s\in[0; 1]\\
t\ne s\end{subarray}}
\Big|\ln\frac{\lambda(t)-\lambda(s)}{t-s}\Big|\leq\ve\},
\end{align*}
$\Lambda$ being a space of strictly increasing continuous mapping from $[0; 1]$ onto itself. 
In the space $(\cD([0; 1]))^N$ the Borel $\sigma$-field is considered. 
\begin{proposition}
\label{prop4'}
For every $u$ the sequence $\{X^{(n)}(u)\}_{n\geq1}$ is weakly compact in $\cD([0; 1]).$
\end{proposition}
Proof.
For a proof we use a test of weak compactness in $\cD([0; 1])$~\cite{billingsley}. 
To do that, one may estimate the difference of $m^{(n)}(u)-X^{(n)}(u)$ by using Propositions \ref{prop2'} and \ref{prop:characteristic_of_diff_m}.
\qed

\begin{corollary}
\label{corl1}
For any $u_1, \ldots, u_N$ a sequence $\{(X^{(n)}(u_1), \ldots,X^{(n)}(u_N) )\}_{n\geq1}$ is weakly compact in $(\cD([0; 1]))^N.$
\end{corollary}

For $u_1< \ldots<u_N,$ we denote an arbitrary weak limit point of $\{(X^{(n)}(u_1), \ldots,X^{(n)}(u_N) )\}_{n\geq1}$ by $(X(u_1), \ldots,X(u_N)).$  We will prove that
$$
(X(u_1), \ldots,X(u_N))\overset{d}{=}(Y^a(u_1), \ldots, Y^a(u_N))
$$
in $(\cD([0; 1]))^N,$ which is our main result, Theorem \ref{thm1}. 

\begin{remark}
 \label{generality}
Without loss of generality and in order to simplify the notation we suppose that the sequence $(X^{(n)}(u_1), \ldots,X^{(n)}(u_N))\Rightarrow(X(u_1), \ldots,X(u_N)),$ $n\rightarrow \infty$. 
\end{remark}

Propositions to follow are checking that the process $(X(u_1), \ldots,X(u_N))$ satisfies all conditions of Definition \ref{defn2}. 

Suppose $B$ is a standard Brownian motion started at 0. For fixed $u$ define processes $\{y^{(n)}(u)\}_{n\geq1}$ in the following way. Consider the partitions $\{t^{(n)}_0, \ldots, t^{(n)}_{N^{(n)}}\},$ $n\in\mbN,$ and put for $t\in[t^{(n)}_k; t^{(n)}_{k+1}), k\in\ov{0;N^{(n)}-1},$
\begin{equation}
\label{prelimit.processes.representation}
\begin{cases}
y^{(n)}_t(u)=B_t-B_{t^{(n)}_k}+\mathop{\circ}\limits^k_{j=1}A_{t^{(n)}_j-t^{(n)}_{j-1}}\left(B_{t^{(n)}_j}-B_{t^{(n)}_{j-1}}+\cdot\right)(u),\\
y^{(n)}_1(u)=\lim_{t\nearrow1}y^{(n)}_t(u).\\
\end{cases}
\end{equation}

Here, for any $j,$ the mapping $\cdot + B_{t^{(n)}_j}-B_{t^{(n)}_{j-1}}$ is a random function from $\mbR$ into inself and is, effectively, a random shift in $\mbR$ by a value $B_{t^{(n)}_j}-B_{t^{(n)}_{j-1}}.$ 

Recall that 
$$
m^{(n)}_t(u)=X^{(n)}_t(u)-\sum_{k: t^{(n)}_k\leq t}\Delta^{(n)}_k(u), \ t\in[0; 1], n\in N.
$$

\begin{proposition}
\label{prop5'}
For any $u$ for any $n$  $(X^{(n)}(u), m^{(n)}(u))\stackrel{d}{=}(y^{(n)}(u), B)$ 
in $(\cD([0; 1]))^2,$ and 
$$
E\int^1_0(y^{(n)}_t(u)-y_t(u))^2dt\to0, n\to\infty, \\
$$
where $y(u)$ is a solution to
$$
\begin{cases}
dy_t(u)=a(y_t(u))dt+dB_t, \\
y_0=u.
\end{cases}
$$
\end{proposition}
Proof. The first assertion follows from the definition of the processes $X^{(n)}, m^{(n)}$, ${n\ge1,}$ while the second one is proved in~\cite{kotelenez}[Propositions 2.5-2.8].
\qed
\begin{corollary}
 \label{diffusion.representation}
 For any $i$
$X(u_i)$ is a diffusion with a unit diffusion coefficient  and a drift coefficient $a,$ that is, for some Brownian motion $B(u_i),$
\begin{equation}
\label{limit.process.is.diffusion}
X_t(u_i)=u_i+\int^t_0a(X_s(u_i))ds+B_t(u_i), \ t\in[0; 1].
\end{equation}
\end{corollary}
Proof. It is left to show that $y^{(n)}\Rightarrow y$ in $\cD([0; 1]),$ as then $X(u_i)\stackrel{d}{=}y$ in $\cD([0; 1]).$ The arguments are standard and thus omitted. 
\qed
\begin{corollary}
\label{corl2}
$$
P\{(X(u_1), \ldots, X(u_N))\in (C([0; 1]))^N\}=1.
$$
\end{corollary}

Now we are to characterize martingale components of $(X(u_1), \ldots, X(u_N))$. For that, define for $f = (f^{1}, f^{2})\in\cD([0;1])$
$$
\theta(f) \equiv \theta(f^1, f^2) = \inf\{1; s\mid f^{1}_s - f^{2}_s \ge 0\}.
$$

We want to prove that for any $i, j$ the process $\big(B_t(u_i), B_t(u_j)\big),$ with $u_i < u_j,$ is a martingale with respect to its own filtration, 
and
$$
\langle B_t(u_i), B_t(u_j)\rangle_{t} = \big( t - \theta(X(u_i), X(u_j)))\big)_{+} = \int_0^t \1_{\{X_s (u_i) = X_s(u_j)\}} ds.
$$
While the first assertion can be verified via standard arguments, the latter needs a more refined approach since we have to check if the processes $X(u_i)$ and  $X(u_j)$ stay equal after $\theta(X(u_i), X(u_j)).$ One technical difficulty is that
$$
\left\{(f^1, f^2)\mid f^1(\theta(f^1, f^2)+\cdot)=f^2(\theta(f^1, f^2)+\cdot)\right\}
$$
is not a closed set in $(\cD([0; 1]))^2$ w.r.t topology induced by the Skorokhod distance. Therefore we cannot use a standard reasoning based on the Portmanteau theorem. Next statements help to overcome this difficulty.

\begin{remark}
In~\cite{konarovskyi}[Lemmas 2.10-2.13] a similar problem for a system of coalescing processes is approached via different arguments based mainly on the theory of martingales. However, our situation differs from that treated in~\cite{konarovskyi} due to the presence of an incorporated jump process and the absence of the convergence of processes stopped on hitting the origin, and while the first problem can be overcomed within the same martinale technique, the second one requires additional reasoning, which we replace with direct calculations of the next proposition and follow-up ones.
\end{remark}

\begin{proposition}
\label{prop7'}
The next event has probability $0:$
$$
\exists i\in\{1, \ldots, N-1\} \exists t\in(0; 1] X_t(u_i)=X_t(u_{i+1}) \mbox{\ and\ } \!\!\sup_{s\in[t; 1]}(X_s(u_{i+1})-X_s(u_i))>0.
$$
\end{proposition}
Proof.
Let $\cD^+([0; 1])=\cD([0; 1])\cap\{f|\inf_{r\in[0; 1]}f_r\geq0\}.$
Define
\begin{align*}
& \Gamma^\delta_\ve=\left\{f\in\cD^+([0; 1])|\exists t\in[0; 1]  \ f(t)<\ve, \int^1_tf(r)dr>\delta\right\},
\\
& \Gamma^\delta=\left\{f\in\cD^+([0; 1])|\exists t\in[0; 1]  \ f(t)=0, \int^1_tf(r)dr>\delta\right\}.
\end{align*}
For $i=\ov{1;N-1},$ denote by $P^{(n), i}$ the distribution of $\Delta X^{(n), i}=X^{(n)}(u_{i+1})-X^{(n)}(u_{i});$ by $P^{(\infty), i},$ the distribution of $\Delta X^i=X(u_{i+1})-X(u_i).$
It is enough to prove that for every fixed $i$ and $\delta>0 \ P^{(\infty),i}(\Gamma^{\delta})=0.$ Since $\Gamma^\delta\subset\cap_{\ve\geq0}\Gamma^\delta_{\ve},$ and $\Gamma_{\ve_1}\subset \Gamma_{\ve_2},\, \ve_1<\ve_2$ we have
\begin{equation*}
P^{(\infty),i}(\Gamma^{\delta})\leq \varliminf_{\ve\to0+}P^{(\infty),i}(\Gamma^{\delta}_\ve).
\end{equation*}

Note that the set $\Gamma^\delta_\ve$ is open in $\cD([0; 1]).$ The Portmanteau theorem~\cite{billingsley} implies  that 
\begin{equation*}
P^{(\infty),i}(\Gamma^\delta_\ve)\leq\varliminf_{n\to\infty}P^{(n),i}(\Gamma^\delta_\ve),
\end{equation*}
since 
$$
(X^{n}(u_{i+1}), X^{(n)}(u_i))\Rightarrow(X(u_{i+1}), X(u_i)), n\to\infty,
$$
and $\Delta X^{(n),i}\Rightarrow\Delta X^i, n\to\infty$ in $\cD([0; 1])$~\cite{billingsley}. Hence
\begin{equation}
 \label{101}
 P^{(\infty),i}(\Gamma^{\delta})\leq \varliminf_{\ve\to0+}\varliminf_{n\to\infty}P^{(n),i}(\Gamma^\delta_\ve). 
\end{equation}

We need an estimate on the modulus of continuity of $\Delta X^{(n),i}.$ For a function $f$ define its modulus of continuity $\omega(f, a) = \sup_{s,t,\ |s-t|\le a} |f_s - f_t|.$ For $s,t\in[t^{(n)}_k;t^{(n)}_{k+1})$
\begin{align*}
|\Delta X^{(n),i}_s-\Delta X^{(n),i}_t| & \leq |X^{(n)}_s(u_i)-X^{(n)}_t(u_i)|+
|X^{(n)}_s(u_{i+1})-X^{(n)}_t(u_{i+1})|
\\
& = |m^{(n)}_s(u_i)-m^{(n)}_t(u_i)|+
|m^{(n)}_s(u_{i+1})-m^{(n)}_t(u_{i+1})|
\\
& \le \omega(m^{(n)}(u_i),\lambda^{(n)})+\omega(m^{(n)}(u_{i+1}),\lambda^{(n)}),
\end{align*}
since $X^{(n)}_s(u_i)-X^{(n)}_t(u_i)=m^{(n)}_s(u_i)-m^{(n)}_t(u_i)$ and the analogous holds for $X^{(n)}(u_{i+1})$.
Define 
$$
\omega^{(n)}=\omega(m^{(n)}(u_i),\lambda^{(n)})+\omega(m^{(n)}(u_{i+1}),\lambda^{(n)}).
$$ 
All $m^{(n)}$ are Brownian motions by Proposition~\ref{prop:characteristic_of_diff_m}, hence $\omega^{(n)}\rightarrow 0, n\rightarrow \infty$ in probability.

If $\inf_{r\in[t^{(n)}_{k-1}; t^{(n)}_k)}\Delta X^{(n),i}_r<\ve$ then
\begin{align*}
\sup_{r\in[t^{(n)}_{k-1}; t^{(n)}_k)}\Delta X^{(n),i}_r & \leq\ve+\omega^{(n)} ,
\\
\Delta X^{(n),i}_{t^{(n)}_k-}=\lim_{r\nearrow t^{(n)}_k}\Delta X^{(n),i}_r & \leq\ve+\omega^{(n)}.
\end{align*}

If for some $t\in[t^{(n)}_{k-1}; t^{(n)}_k)$\ $\Delta X^{(n),i}_t<\ve $ and
$\int^1_t\Delta X^{(n),i}_rdr>\delta,$ then 
\begin{align*}
\int^{t^{(n)}_k}_t\Delta X^{(n),i}_rdr & \leq (t^{(n)}_k-t)\sup_{r\in[t^{(n)}_{k-1}; t^{(n)}_k)}\Delta X^{(n),i}_r\leq \lambda^{(n)} (\ve+\omega^{(n)})\leq \ve+\omega^{(n)},
\\
\int_{t^{(n)}_k}^1\Delta X^{(n),i}_rdr & > \delta-\int^{t^{(n)}_k}_t\Delta X^{(n),i}_rdr\geq\delta-\lambda^{(n)} (\ve+\omega^{(n)})\geq \delta - (\ve+\omega^{(n)}).
\end{align*}
Using these two remarks:
\begin{align}
\label{eq11'}
&P^{(n),i}(\Gamma^\delta_\ve)  =P\left\{\exists t: \Delta X^{(n),i}_t<\ve, \int^1_t\Delta X^{(n),i}_rdr>\delta\right\}  \nonumber
\\
& = 
\sum^{N^{(n)}}_{k=1}
P\left\{\inf_{r\in[0; t^{(n)}_{k-1})}\Delta X^{(n),i}_r\geq\ve; \exists t\in[t^{(n)}_{k-1}; t^{(n)}_k) \!\Delta X^{(n),i}_t<\ve, \int_{t}^1\Delta X^{(n),i}_rdr>\delta\right\} \nonumber
\\
& = \sum^{N^{(n)}}_{k=1}
P\left\{\inf_{r\in[0; t^{(n)}_{k-1})}\Delta X^{(n),i}_r\geq\ve; \inf_{r\in[t^{(n)}_{k-1}; t^{(n)}_{k})}\Delta X^{(n),i}_r<\ve; \int_{t}^1\Delta X^{(n),i}_rdr>\delta;\right. \nonumber
\\
&\left. \phantom{abcabcabc} \ve+\omega^{(n)}\geq 2\ve\wedge \frac{\delta}{2}\right\} \nonumber \\  
& \phantom{abc} + 
\sum^{N^{(n)}}_{k=1}
P\left\{\inf_{r\in[0; t^{(n)}_{k-1})}\Delta X^{(n),i}_r\geq\ve; \inf_{r\in[t^{(n)}_{k-1}; t^{(n)}_{k})}\Delta X^{(n),i}_r<\ve;
\int_{t}^1\Delta X^{(n),i}_rdr>\delta; \right.
\nonumber
\\ 
& \phantom{abcabcabc} \left.
\ve+\omega^{(n)}<2\ve\wedge \frac{\delta}{2}\right\} \nonumber
\\
&\leq 
P\left\{\ve+\omega^{(n)}\geq 2\ve\wedge \frac{\delta}{2}\right\}+ \sum^{N^{(n)}}_{k=1} 
P\left\{\inf_{r\in[0; t^{(n)}_{k-1})}\Delta X^{(n),i}_r\geq\ve;\right. 
 \nonumber
\\ 
& \phantom{abcabcabc} \inf_{r\in[t^{(n)}_{k-1}; t^{(n)}_{k})}\Delta X^{(n),i}_r<\ve;
\int_{t}^1\Delta X^{(n),i}_rdr>\delta; \Delta X^{(n),i}_{t^{(n)}_k-}\leq2\ve; \nonumber \\
& \phantom{abcabcabc}  \left.
\int_{t^{(n)}_k}^1\Delta X^{(n),i}_rdr\geq\delta/2; 
\ve+\omega^{(n)}<2\ve\wedge \frac{\delta}{2}\right\} \nonumber
\\ 
& \leq
P\left\{\ve+\omega^{(n)}\geq 2\ve\wedge \frac{\delta}{2}\right\}+
\sum^{N^{(n)}}_{k=1} 
P\left\{\!\inf_{r\in[0; t^{(n)}_{k-1})}\!\Delta X^{(n),i}_r\geq\ve; 
 \right. \nonumber
\\ 
& \phantom{abcabcabc} \left. \inf_{r\in[t^{(n)}_{k-1}; t^{(n)}_{k})}\!\Delta X^{(n),i}_r<\ve;
\Delta X^{(n),i}_{t^{(n)}_k-}\leq2\ve;
\int_{t^{(n)}_k}^1\Delta X^{(n),i}_rdr\geq\delta/2\right\} \nonumber
\\
& =
P\left\{\ve+\omega^{(n)}\geq 2\ve\wedge \frac{\delta}{2}\right\}
+ 
E\sum^{N^{(n)}}_{k=1} 
P\left(
\inf_{r\in[0; t^{(n)}_{k-1})}\Delta X^{(n),i}_r\geq\ve; \right.
 \nonumber
\\
& \phantom{abcabcababc} \left. \Delta X^{(n),i}_{t^{(n)}_k-}\leq2\ve; \Delta X^{(n),i}_{t^{(n)}_k-}\leq2\ve;
\int_{t^{(n)}_k}^1\Delta X^{(n),i}_rdr\geq\delta/2
\left\vert\cF_{t^{(n)}_k}\right.\right)
\nonumber
\\
& =
P\left\{\ve+\omega^{(n)}\geq 2\ve\wedge \frac{\delta}{2}\right\}
+ 
E\sum^{N^{(n)}}_{k=1} 
\1_{\{\inf_{r\in [0; t^{(n)}_{k-1})}\Delta X^{(n),i}_r\geq\ve\}}
\nonumber
\\
& \phantom{abc} \cdot 
\1_{\{\inf_{r\in [t^{(n)}_{k-1}; t^{(n)}_{k})}\Delta X^{(n),i}_r<\ve\}}
\1_{\{\Delta X^{(n),i}_{t^{(n)}_k-}\leq2\ve \}}
P\left( \int_{t^{(n)}_k}^1\Delta X^{(n),i}_rdr\geq\delta/2\left\vert\cF_{t^{(n)}_k}\right.\right).
\end{align}
To proceed further, consider, for $t\in[t^{(n)}_j; t^{(n)}_{j+1}),$
\begin{align*}
E\left(\Delta X^{(n),i}_t\left\vert\cF_{t^{(n)}_j}\right.\right) & =
E\left(\vf_{t^{(n)}_{j}, t}
(X^{(n)}_{t^{(n)}_{j}}(u_{i+1}))
-
\vf_{t^{(n)}_{j}, t}
(X^{(n)}_{t^{(n)}_{j}}(u_i))\left\vert\cF_{t^{(n)}_j}\right.\right)
\\
& =E\left(\vf_{t^{(n)}_{j}, t}
(X^{(n)}_{t^{(n)}_{j}}(u_{i+1}))-\vf_{t^{(n)}_{j}, t}
(X^{(n)}_{t^{(n)}_{j}}(u_i))\left\vert\cF_{t^{(n)}_j}\right.\right)
\\
& = 
E\left(\vf_{t^{(n)}_{j}, t}
(\xi_1)-\vf_{t^{(n)}_{j}, t}
(\xi_2)\right)\Big|_{\xi_1=X^{(n)}_{t^{(n)}_{j}}(u_{i+1}), \xi_2=X^{(n)}_{t^{(n)}_{j}}(u_i)}.
\end{align*}
Here, from the definition of the Brownian web, 
$$
E\left(\vf_{t^{(n)}_{j}, t}
(\xi_1)-\vf_{t^{(n)}_{j}, t}
(\xi_2)\right)=E\sqrt{2}B_{(t-t^{(n)}_j)\wedge\theta}(\frac{\xi_2-\xi_1}{\sqrt{2}}),
$$
where $\{B_s(x)|s\geq0\}, x\in\mbR^+,$ is a Brownian motion started at $x,$ and $\theta$ is its moment of hitting  $0$. Thus
\begin{equation*}
\label{eq12'}
E\left(\Delta X^{(n),i}_t\left\vert\cF_{t^{(n)}_j}\right.\right)\leq\Delta X^{(n),i}_{t^{(n)}_j}.
\end{equation*}

The Gronwall--Belmann lemma implies that
\begin{equation*}
\label{eq13'}
\Delta X^{(n),i}_{t^{(n)}_j}\leq \Delta X^{(n),i}_{t^{(n)}_j-}\cdot \e^{C_a(t^{(n)}_{j+1}-t^{(n)}_j)}.
\end{equation*}

Since 
\begin{align*}
 \Delta X^{(n),i}_{t^{(n)}_{j}-} & =\lim_{r\nearrow t^{(n)}_{j}}\Delta X^{(n),i}_{r}=\lim_{r\nearrow t^{(n)}_{j}} \left[\vf_{t^{(n)}_{j-1}, r} (X^{(n)}_{t^{(n)}_{j-1}}(u_{i+1}))-\vf_{t^{(n)}_{j-1}, r}(X^{(n)}_{t^{(n)}_{j-1}}(u_i))\right]
\\
& =\vf_{t^{(n)}_{j-1}, t^{(n)}_{j}} \left(X^{(n)}_{t^{(n)}_{j-1}}(u_{i+1}))-\vf_{t^{(n)}_{j-1}, t^{(n)}_{j}}(X^{(n)}_{t^{(n)}_{j-1}}(u_i)\right)
\end{align*}
we end with
\begin{equation*}
\label{eq14'} 
E\left( \Delta X^{(n),i}_{t} \left\vert\cF_{t^{(n)}_k}\right.\right)\leq e^{C_a(1-t^{(n)}_k) } \Delta X^{(n),i}_{t^{(n)}_{k}-} 
\leq  e^{C_a } \Delta X^{(n),i}_{t^{(n)}_{k}-}.
\end{equation*}
This gives: 
\begin{align*}
& P\left( \int_{t^{(n)}_k}^1\Delta X^{(n),i}_rdr\geq\delta/2 \left\vert\cF_{t^{(n)}_k} \right.\right)\leq 
\frac{2}{\delta}E\left(\int_{t^{(n)}_k}^1\Delta X^{(n),i}_rdr\left\vert\cF_{t^{(n)}_k}\right.\right)
\\
&\leq
\frac{2}{\delta}
E\left(
\varliminf_{n\to\infty} 
\sum^{m-1}_{j=0}\Delta X^{(n),i}_{(1-t^{(n)}_k)\frac{j}{m}+t^{(n)}_k}\frac{1-t^{(n)}_k}{m}\left\vert\cF_{t^{(n)}_k}\right.\right)
\\
&\leq \frac{2}{\delta}
\varliminf_{n\to\infty} 
\sum^{m-1}_{j=0} E\left(\Delta X^{(n),i}_{(1-t^{(n)}_k)\frac{j}{m}+t^{(n)}_k}\frac{1-t^{(n)}_k}{m}\left\vert\cF_{t^{(n)}_k}\right.\right)
\\
&\leq
\frac{2}{\delta}
\varliminf_{n\to\infty}
\sum^{m-1}_{j=0}\e^{C_a}\Delta X^{(n),i}_{t^{(n)}_k-}\frac{1-t^{(n)}_k}{m}\leq
\frac{2}{\delta}\e^{C_a}\Delta X^{(n),i}_{t^{(n)}_k-}.
\end{align*}
Substituting the last estimate into \eqref{eq11'}, we get:
\begin{align*}
P^{(n),i}(\Gamma^\delta_\ve) & \leq 
P\left\{\ve+\omega^{(n)}\geq 2\ve\wedge \frac{\delta}{2}\right\}
+
E\sum^{N^{(n)}}_{k=1} 
\1_{\{\inf_{r\in[0; t^{(n)}_{k-1})}\Delta X^{(n),i}_r\geq\ve\}}
\\ & \phantom{abcabc} \cdot
\1_{\{\inf_{r\in[t^{(n)}_{k-1}; t^{(n)}_{k})}\Delta X^{(n),i}_r<\ve\}}
\1_{\{\Delta X^{(n),i}_{t^{(n)}_k-}\leq2\ve \}} \cdot
\frac{2}{\delta}\e^{C_a}\Delta X^{(n),i}_{t^{(n)}_k-}
\\
& \leq
P\left\{\ve+\omega^{(n)}\geq 2\ve\wedge \frac{\delta}{2}\right\}
\\ & \phantom{abc} +  
\sum^{N^{(n)}}_{k=1} 
P\left\{\inf_{r\in[0; t^{(n)}_{k-1})}\Delta X^{(n),i}_r\geq\ve; \inf_{r\in[t^{(n)}_{k-1}; t^{(n)}_{k})}\Delta X^{(n),i}_r<\ve\right\}\cdot
\frac{4}{\delta}\e^{C_a} \ve
\\  & =
P\left\{\ve+\omega^{(n)}\geq 2\ve\wedge \frac{\delta}{2}\right\}+
\frac{4}{\delta}\e^{C_a}\cdot \ve.
\end{align*}

Thus in (\ref{101}) we have: 
$$
\lim_{\ve\to0+}\varliminf_{n\to\infty}P^{(n),i}(\Gamma^\delta_\ve)=0.
$$
\qed

\begin{proposition}
 \label{hitting.moment.continuity}
 Let $a^{(1)}, a^{(2)}$ be Lipschitz continuous functions on the real line. Let $\xi^{(1)}, \xi^{(2)}$ be solutions on $[0; 1]$ to the following SDEs:
$$
\begin{cases}
d\xi^{(k)}_t= w^{(k)}_t+a^{(k)}(\xi^{(k)}_t)dt, \ \xi^{(k)}_0=u^{(k)}, \\
k=1,2,  \ u^{(1)}\leq u^{(2)}, 
\end{cases}
$$
where $w^{(1)}, w^{(2)}$ are independent standard Brownian motion started at 0. 
Then the function $\theta$ is $\mathrm{Law}(\xi^{(1)}, \xi^{(2)})-$a.s. continuous on $C([0;1], \mbR^2)$.
\end{proposition}
Proof.
Since the reasoning is standard we give only a sketch of the proof. Define $\xi = (\xi^{(1)}, \xi^{(2)}).$  Suppose $g^{(n)}\rightarrow\xi$ in $(\cD([0;1]))^2, n\rightarrow\infty.$ Put $\theta^{(n)} = \theta(g^{(n)}), \theta^{\infty} = \theta(\xi).$ Then one can prove that $\xi^{(1)}(\varliminf_{n\to\infty}\theta^{(n)}) -\xi^{(2)}(\varliminf_{n\to\infty}\theta^{(n)}) = 0,$ so 
 $$
 \theta^{\infty}\le \varliminf\limits_{n\rightarrow\infty}\theta^{(n)}.
 $$
The moment $\theta^{\infty}$ is easily seen to be a Markov moment with respect to the filtration generated by $\xi$~\cite{kallenberg} so for the process $\xi^{(2)}_{\theta^{\infty}+\cdot}-\xi^{(1)}_{\theta^{\infty}+\cdot}$ we have
$$
d\left(\xi^{(2)}_{\theta^{\infty}+t}-\xi^{(1)}_{\theta^{\infty}+t}\right)= \sqrt{2} b_t+\left(a^{(2)}(\xi^{(2)}_{\theta^{\infty}+t})-a^{(1)}(\xi^{(1)}_{\theta^{\infty}+t})\right)dt,
$$ 
where $b$ is a standard Brownian motion. The Girsanov theorem~\cite{random.proc.stats} implies that $\xi^{(2)}_{\theta^{\infty}+\cdot}-\xi^{(1)}_{\theta^{\infty}+\cdot}$ obeys the iterated logarithm law. So  
 $$
 \varliminf_{h\searrow 0}\frac{\xi^{(2)}_{\theta^{\infty} + h} - \xi^{(1)}_{\theta^{\infty} + h}}{\sqrt{2h\log\log h^{-1}}} = 1 \mbox{\ a.s.}.
 $$
 Thus
 $$
 \theta^{\infty}\ge \varlimsup\limits_{n\rightarrow\infty}\theta^{(n)}.
 $$
\qed

Recall
$$
\tau^{(n)} = \inf\left\{1; s\mid X^{(n)}_s(u_1) =  X^{(n)}_s(u_2)\right\}, n\ge1,
$$
$u_1 < u_2,$ and define $\tau = \inf\left\{1; s\mid X_s(u_1) =  X_s(u_2)\right\}.$
\begin{proposition}
 \label{cor:weak.conv.hitting.times}
$\tau^{(n)}\Rightarrow\tau, n\to\infty.$ 
\end{proposition} 
Proof. By the Skorokhod representation theorem~\cite{kallenberg} we may assume that with probability $1$
$$
(X^{(n)}(u_1), X^{(n)}(u_2)) \rightarrow (X(u_1), X(u_2)), n\to\infty,
$$
in  $(\cD([0; 1]))^2,$ for some  $(X(u_1), X(u_2))$ described in Corollary~\ref{diffusion.representation}. In this notation, we need to establish the following:
$$
\theta(X^{(n)}(u_1), X^{(n)}(u_2))\Rightarrow \theta(X(u_1), X(u_2)), n\to\infty.
$$ 
Since both $X(u_1)$ and $X(u_2)$ by Corollary~\ref{corl2} are continuous a.s., $(X^{(n)}(u_1),$ $X^{(n)}(u_2)$ converge to $(X(u_1), X(u_2))$ uniformly. So it easy to see that 
\begin{equation}
\label{star}
\theta(X(u_1), X(u_2)) \le \varliminf\limits_{n\rightarrow\infty}\theta(X^{(n)}(u_1), X^{(n)}(u_2)).
\end{equation}

Either coordinate of $(X^{(n)}(u_1), X^{(n)}(u_2))$ admits the representation of~(\ref{prelimit.processes.representation}), with $B$ replaced with some Brownian motions $m^{(n)}(u_k), k = 1,2,$ respectively. 

Let $\tilde{B}$ be a standard Brownian motion  with $B_0=0,$ independent of all the processes used before, the initial probability space being extended in a usual way, if needed. 
Define 
\begin{align*}
\tilde{m}^{(n)}_t(u_2) & = m^{(n)}_t(u_2) \1_{\{t<\tau^{(n)}\}} + \left( \tilde{B}_t - \tilde{B}_{\tau^{(n)}} + m^{(n)}_{\tau^{(n)}}(u_2) \right) \1_{\{t\ge\tau^{(n)}\}} \\
& = m^{(n)}_{t\wedge\tau^{(n)}}(u_2) + \tilde{B}_t - \tilde{B}_{t\wedge \tau^{(n)}},
\end{align*}
and consider the process $\tilde{X}^{(n)}(u_2)$ constructed as in Equation~\ref{prelimit.processes.representation} with $\tilde{m}^{(n)}(u_2)$ substituted for $B.$ Then one can easily check that
\begin{equation}
\label{eq:rebuilt.prelimit.processses} 
\tilde{X}^{(n)}_t(u_2) \1_{\{t<\tau^{(n)}\}} = X^{(n)}_t(u_2) \1_{\{t<\tau^{(n)}\}}.
\end{equation}
Consider new Brownian motions $b(u_1)$ and $b(u_2),$ independent of all processes used above. Construct processes $y^{(n)}(u_1), y^{(n)}(u_2)$ using the representation of~\ref{prelimit.processes.representation} with $b(u_1)$ or $b(u_2)$ instead of $B,$ respectively. Then
\begin{equation}
\label{independent.times.equals.ours}
(X^{(n)}(u_1), \tilde{X}^{(n)}(u_2)) \stackrel{d}{=} (y^{(n)}(u_1), y^{(n)}(u_2)).
\end{equation}
By Proposition~\ref{prop5'} there exist processes $(y(u_1), y(u_2))$ such that
$$
\begin{cases}
d{y}_t(u_k)=a({y}_t(u_k))dt+d\tilde{m}_t(u_k), \\
{y}_0(u_k)=u_k, k  = 1,2,
\end{cases}
$$
and, for each $k,$ the sequence $\left\{ y^{(n)}(u_k) \right\}_{n\ge1} $ converges to $y(u_k)$ in $\cL_{2}(\Omega\times(0;1)),$ and, consequently, in $\cD([0; 1]).$  Here $\tilde{m}(u_1)$ and $\tilde{m}(u_2)$ are standard Brownian motions started from $0.$ Since prelimit processes $b(u_1)$ and $b(u_2)$ are independent one can prove 
$\tilde{m}(u_1)$ and $\tilde{m}(u_2)$ to be independent, too. So such are the processes $y(u_1)$ and $y(u_2).$ 

Thus we have, by~(\ref{eq:rebuilt.prelimit.processses}) and~(\ref{independent.times.equals.ours}):
\begin{align*}
\tau^{n} & \stackrel{d}{=} \theta(X^{(n)}(u_1), X^{(n)}(u_2)) = \theta(X^{(n)}(u_1), \tilde{X}^{(n)}(u_2))  \\
& \stackrel{d}{=} \theta(y^{(n)}(u_1), y^{(n)}(u_2)), n\in \mbN.
\end{align*}
By Proposition~\ref{hitting.moment.continuity} the mapping $\theta$ is a.e. continuous w.r.t the distribution of $(y(u_1), y(u_2)),$ which implies that
$$
 \theta(y^{(n)}(u_1), y^{(n)}(u_2)) \rightarrow \theta(y(u_1), y(u_2)), n\to\infty, \mbox{\ a.s.},
$$
so $\tau^{n}\Rightarrow\theta(y(u_1), y(u_2)), n\to\infty.$ It is to left to prove that
$$
\theta(X(u_1), X(u_2)) \stackrel{d}{=} \theta(y(u_1), y(u_2)).
$$
For that recall that, by Corollary~\ref{prelimit.processes.representation},
$$
\begin{cases}
\label{starstar}
d{X}_t(u_k)=a({X}_t(u_k))dt+d m_t(u_k), \\
{X}_0(u_k)=u_k, k  = 1,2,
\end{cases}
$$
where $m(u_1)$ and $m(u_2)$ are Brownian motions. For any $t,$
$$
\int_{0}^t a(X^{(n)}_s(u_k))ds \rightarrow \int_{0}^t a(X_s(u_k))ds,
$$
so one can prove that $m^{(n)}(u_k)\rightarrow m(u_k), n\rightarrow\infty,$ uniformly, for each $k.$ Then Proposition~\ref{prop:characteristic_of_diff_m} implies that 
\begin{align*}
\left\langle m(u_1), m(u_2) \right\rangle_t & = \lim_{n\rightarrow\infty} \left\langle m^{(n)}(u_1), m^{(n)}(u_2) \right\rangle_t 
\\
& =
\lim_{n\rightarrow\infty} \left(t - \theta(X^{(n)}(u_1), X^{(n)}(u_2)\right)_{+} \\
& \le  \left(t - \varlimsup_{n\rightarrow\infty}\theta(X^{(n)}(u_1), X^{(n)}(u_2)\right)_{+}.
\end{align*}
Using~\ref{star}, we have that
$$
0 \le \left\langle m(u_1), m(u_2) \right\rangle_t \le \left(t - \theta(X(u_1), X(u_2)\right)_{+},
$$
so 
$$
\left\langle m_{\cdot \wedge \theta(X(u_1), X(u_2)}(u_1), m_{\cdot \wedge \theta(X(u_1), X(u_2)}(u_2) \right\rangle_t = 0, t\in [0;1].
$$
Thus $\theta(X(u_1), X(u_2)) \stackrel{d}{=} \theta(y(u_1), y(u_2)).$ Indeed, either $X(u_k)$ is a strong unique solution of Equation~\ref{starstar}, so the process $X_{\cdot \wedge \theta(X(u_1), X(u_2)}(u_k)$ is adapted to the filtration generated by  $m_{\cdot \wedge \theta(X(u_1), X(u_2)}(u_k).$ Since processes $m_{\cdot \wedge \theta(X(u_1), X(u_2)}(u_1)$ and $m_{\cdot \wedge \theta(X(u_1), X(u_2)}(u_2)$ are independent this finishes the proof.
\qed

Recall from Equation~\ref{limit.process.is.diffusion}
$$
X_t(u_i)=u_i+\int^t_0a(X_s(u_i))ds+B_t, \ t\in[0; 1],
$$
for any $i.$ We will write $B_t(u_i)$ to emphasize that the Brownian motion in Equation~\ref{limit.process.is.diffusion} depends on the starting point $u_i.$
\begin{proposition}
 \label{limit.characteristic}
For any $i, j\colon u_i<u_j,$ the process $(B_t(u_i), B_t(u_j))$ is a martingale with respect to its own filtration,
and
$$
\left\langle B_t(u_i), B_t(u_j)\right\rangle_{t} = \left( t - \tau\right)_{+} = \int_0^t \1_{\{X_s (u_i) = X_s(u_j)\}} ds ,
$$
where $\tau = \inf\{1; s\mid X_s(u_1) =  X_s(u_2)\}.$
\end{proposition}
Proof.
A standard approach can be used to check that going to the limit $(B_t(u_i), B_t(u_j))$ preserves  the martingale property of $(m^{(n)}(u_i), m^{(n)}(u_j))$ . Also 
$$
\left\langle B_t(u_i), B_t(u_j)\right\rangle_{t} = \lim_{n\rightarrow\infty} \left\langle m^{(n)}(u_i), m^{(n)}(u_j)\right\rangle_t.
$$
Proposition~\ref{cor:weak.conv.hitting.times} gives, for any $t\in[0;1],$ that
$$
E\big( t - \tau^{(n)}\big)_{+} \rightarrow E\big( t - \tau\big)_{+}, n\to\infty,
$$ 
and Proposition~\ref{prop7'} implies that
$$
\big( t - \tau\big)_{+} = \int_0^t \1_{\{X_s (u_i) = X_s(u_j)\}} ds.
$$
\qed

\section{Main Result}
Now we can formulate and prove the main result.

\begin{theorem}
\label{thm1}
Let $u_1<\ldots<u_N.$ Then
$$
(X^{(n)}(u_1), \ldots, X^{(n)}(u_N))\Rightarrow (Y^a(u_1), \ldots, Y^a(u_N)), n\to\infty,
$$
in $(\cD([0; 1]))^N.$
\end{theorem}
Proof. One need use Definition~\ref{defn2}. Here, Proposition~\ref{limit.characteristic} gives the first condition of the definition, while Corollary~\ref{diffusion.representation} is exactly the second condition.
\qed

\section{Appendix}
\label{appendix}
Proof of Proposition~\ref{tech}.
We give a sketch of the proof only for $s=0$ and $t=1$ to keep the notation simple. Throughout the proof the function $\alpha\in(C([0;1]))^N$ and the ordered set $\mathcal{U}=\{u_n\mid n=\overline{1, N}\}$ are fixed. Denote  $\xi_{n,k}=\vf_{\frac{k}{n}, \frac{k+1}{n}}(0), k = \overline{0, n-1}, n\in \mbR.$ Consider 
$$ 
\cF_{s, t} = \sigma\left(\vf_{r_1, r_2}(v), s\leq r_1 \leq r_2\leq t, v\in\mbR\right).
$$

Given an arbitrary ordered set $\mathcal{A}$ with $M$ elements and an arbitrary function $\beta\in (C([0;1]))^M$ it holds that
\begin{equation}
\label{exp:mean}
E\mathcal{E}^M_{s, t}(\mathcal{A}; \beta)=1,  
\end{equation}
which was proved in~\cite{dorogovtsev:density}. Hereinbelow we denote a size of any finite set $\mathcal{A}$ as $|A|.$ 

For $k=\overline{0,n-1}$ define a set $\mathcal{U}_k=\{u_{k,i}\mid i = \overline{1, M_k}\}$ as a set of all distinct elements of $\{ \vf_{0, \frac{k}{n}}(u_j)\mid j=\overline{1, N}\}$ whose ordering is inherited from that of $\mathcal{U}$ in a obvious way.
Define, for a fixed $k,$
\begin{align*}
l_{k,i} &  = \min\left\{j\mid \vf_{0, \frac{k+1}{n}}(u_j) = u_{k,i}\right\}, i=\overline{1, M_k}, \\
\mathcal{L}_k & = (\mathcal{L}_{k,1}, \ldots, \mathcal{L}_{k,M_k}) = (l_{k,1}, \ldots, l_{k,M_k}).
\end{align*}

Equation~(\ref{eq1}) implies that, for any $k,$
\begin{align*}
\mathcal{E}^N_{0, \frac{k+1}{n}} & (\mathcal{U}; \alpha) = \mathcal{E}^N_{0, \frac{k}{n}}(\mathcal{U}; \alpha) 
\\
& \phantom{abc} \cdot
\mathrm{exp}\left( \sum_{j=1}^N\int\limits_{\frac{k}{n}}^{\tau_j(\mathcal{U})\vee \frac{k}{n}}\!\!\alpha_j(r)d\vf_{\frac{k}{n}, r}(\vf_{0, \frac{k}{n}}(u_j)) -\frac{1}{2} \sum_{j=1}^N\int\limits_{\frac{k}{n}}^{\tau_j(\mathcal{U})\vee \frac{k}{n}}\!\!\alpha_j^2(r)dr \right) 
\\
& = 
\mathcal{E}^N_{0, \frac{k}{n}}(\mathcal{U}; \alpha) \\
& \phantom{abc} \cdot 
\mathrm{exp}\left( \sum_{j=1}^{M_k}\int\limits_{\frac{k}{n}}^{\tau_{l_{k,j}}(\mathcal{U})}\!\!\alpha_{l_{k,j}}(r)d\vf_{\frac{k}{n}, r}(\vf_{0, \frac{k}{n}}(u_{l_{k,j}})) -\frac{1}{2} \sum_{j=1}^{M_k}\int\limits_{\frac{k}{n}}^{\tau_{l_{k,j}}(\mathcal{U})}\!\!\alpha_{l_{k,j}}^2(r)dr \right) 
\\
& =
\mathcal{E}^N_{0, \frac{k}{n}}(\mathcal{U}; \alpha) \\
& \phantom{abc} \cdot 
\mathrm{exp}\left( \sum_{j=1}^{M_k}\int\limits_{\frac{k}{n}}^{\tau_j(\mathcal{U}_k)}\!\!\alpha_j(\mathcal{L}_k)(r)d\vf_{\frac{k}{n}, r}(u_{k, j}) -\frac{1}{2} \sum_{j=1}^{M_k}\int\limits_{\frac{k}{n}}^{\tau_j(\mathcal{U}_k)}\!\!\alpha_j(\mathcal{L}_k)^2(r)dr \right),
\end{align*}
where $\alpha(\mathcal{L}_k)=(\alpha_1(\mathcal{L}_k), \ldots, \alpha_j(\mathcal{L}_k))$ is a random $\cF_{0, \frac{k}{n}}-$measurable function defined by a relation $\alpha_j(\mathcal{L}_k)=\alpha_{l_{k,j}}=\alpha_{\mathcal{L}_{k,j}}, j=\overline{1,M_k}.$ Thus 
\begin{equation*}
\label{exp:decomp}
\mathcal{E}^N_{0, \frac{k+1}{n}}(\mathcal{U}; \alpha) = \mathcal{E}^N_{0, \frac{k}{n}}(\mathcal{U}; \alpha)\cdot  \mathcal{E}^{M_k}_{ \frac{k}{n}, \frac{k+1}{n}}(\mathcal{U}_k; \alpha(\mathcal{L}_k)). 
\end{equation*}
We omit a standard discussion on measurability issues. What is essential is that for a deterministic ordered tuple $\mathcal{L}$ the function $\alpha(\mathcal{L})$ is nonrandom. Since $\mathcal{E}^N_{0, \frac{k}{n}}\left(\mathcal{U}; \alpha\right)$ is $\cF_{0, \frac{k}{n}}-$measurable and for any nonrandom ordered set $\mathcal{A}$ and function $\beta$ a random variable $\mathcal{E}^{|\mathcal{A}|}_{\frac{k}{n}, \frac{k+1}{n}}\left(\mathcal{A}; \beta\right)$ is independent of $\cF_{0, \frac{k}{n}}$ by the construction we have that
\begin{align}
\label{exp:xi.cond}
E\left(\xi_{n,k} \mathcal{E}^N_{0, \frac{k+1}{n}}(\mathcal{U}; \alpha) \left\vert \cF_{0, \frac{k}{n}} \right.\right) & = E\left(\xi_{n,k}\mathcal{E}^N_{0, \frac{k}{n}}(\mathcal{U}; \alpha)\cdot  \mathcal{E}^{|\mathcal{U}_k|}_{ \frac{k}{n}, \frac{k+1}{n}}(\mathcal{U}_k; \alpha(\mathcal{L}_k))  \left\vert \cF_{0, \frac{k}{n}} \right.\right) \nonumber
\\
& =
\mathcal{E}^N_{0, \frac{k}{n}}(\mathcal{U}; \alpha) E\left(\xi_{n,k}\mathcal{E}^{|\mathcal{A}|}_{\frac{k}{n}, \frac{k+1}{n}}(\mathcal{A}; \alpha(\mathcal{K})) \right)\Big|_{\mathcal{A}=\mathcal{U}_k, \mathcal{K}=\mathcal{L}_k}.
\end{align}
Additionally, Equation~(\ref{exp:mean}) yields that
\begin{align}
\label{exp:sol.cond}
E\left(\xi_{n,k} \mathcal{E}^N_{0, 1}(\mathcal{U}; \alpha) \left\vert \cF_{0, \frac{k+1}{n}} \right.\right) & =
\xi_{n,k} \mathcal{E}^N_{0, \frac{k+1}{n}}(\mathcal{U}; \alpha) E\left(\mathcal{E}^{|\mathcal{A}|}_{\frac{k+1}{n}, 1}(\mathcal{A}; \alpha(\mathcal{K}))  \right)\Big|_{\mathcal{A}=\mathcal{U}_k, \mathcal{K}=\mathcal{L}_k} \nonumber
\\ 
& = \xi_{n,k} \mathcal{E}^N_{0, \frac{k+1}{n}}(\mathcal{U}; \alpha).
\end{align}

Using~(\ref{exp:sol.cond}) and (\ref{exp:xi.cond}):
\begin{align}
\label{app:main}
E\sum_{k=1}^{n-1}\xi_{n,k}\mathcal{E}^N_{0, 1}(\mathcal{U}; \alpha) & = \sum_{k=1}^{n-1}EE \left( \xi_{n,k}\mathcal{E}^N_{0, 1}(\mathcal{U}; \alpha) \left\vert \cF_{0, \frac{k+1}{n}}\right.\right) \nonumber 
\\ 
& = 
\sum_{k=1}^{n-1}E \xi_{n,k} \mathcal{E}^N_{0, \frac{k+1}{n}}(\mathcal{U}; \alpha)  = 
\sum_{k=1}^{n-1}E E \left( \xi_{n,k} \mathcal{E}^N_{0, \frac{k+1}{n}}(\mathcal{U}; \alpha) \left\vert \cF_{0, \frac{k}{n}}\right.\right) \nonumber
\\ 
& = \sum_{k=1}^{n-1}E \mathcal{E}^N_{0, \frac{k}{n}}(\mathcal{U}; \alpha) E\left(\xi_{n,k}\mathcal{E}^{|\mathcal{A}|}_{\frac{k}{n}, \frac{k+1}{n}}(\mathcal{A}; \alpha(\mathcal{K}))  \right)\Big|_{\mathcal{A}=\mathcal{U}_k, \mathcal{K}=\mathcal{L}_k}.
\end{align}
We are going to estimate $E\xi_{n,k}\mathcal{E}^{|\mathcal{A}|}_{\frac{k}{n}, \frac{k+1}{n}}(\mathcal{A}; \alpha(\mathcal{K}))$ for fixed ordered sets $\mathcal{A}$ and $\mathcal{K}, \mathcal{K}\subset\{1,\ldots,N\}.$ To start, the stationarity of the Brownian web yields that
\begin{align*}
E\left(\xi_{n,k}\mathcal{E}^{|\mathcal{A}|}_{\frac{k}{n}, \frac{k+1}{n}}(\mathcal{A}; \alpha(\mathcal{K}))  \right) & = 
E \vf_{\frac{k}{n}, \frac{k+1}{n}}(0) \mathcal{E}^{|\mathcal{A}|}_{\frac{k}{n}, \frac{k+1}{n}}(\mathcal{A}; \alpha(\mathcal{K})) 
\\ & =
E \vf_{0, \frac{1}{n}}(0) \mathcal{E}^{|\mathcal{A}|}_{0, \frac{1}{n}}(\mathcal{A}; \alpha(\mathcal{K})).
\end{align*}
Let $\mathcal{A}=\{a_j\mid j=\overline{1,M}\},$ $\alpha(\mathcal{K})=(\beta_1,\ldots,\beta_M), M \le N.$ Note that 
$$
\alpha^{\star} = \sup_{x\in[\min\mathcal{U};\max\mathcal{U}], s\in[0;1]} |\alpha(x, s)| \ge \sup_{s\in[0;1], j=\overline{1,M}} |\beta_j(s)|.
$$ An application of the Ito formula gives, for $t\in[0;\frac{1}{n}]:$
\begin{align*}
 \vf_{0, t}(0) \mathcal{E}^{|\mathcal{A}|}_{0, t}(\mathcal{A}; \beta) = & \int\limits_0^t \! \mathcal{E}^{|\mathcal{A}|}_{0, s}(\mathcal{A}; \beta) d\vf_{0, s}(0) 
\\ 
& + \sum_{j=1}^{M}\int\limits_{0}^{t\wedge \tau_j(\mathcal{A})}  \!\vf_{0, s}(0) \mathcal{E}^{|\mathcal{A}|}_{0, s}(\mathcal{A}; \beta) \beta_j(s)d\vf_{0, s}(a_j) 
\\
& + \sum_{j=1}^{M} \int\limits_{0}^{t\wedge \tau_j(\mathcal{A})} \mathcal{E}^{|\mathcal{A}|}_{0, s}(\mathcal{A}; \beta) d\left\langle \vf_{0, \cdot}(0), \int\limits_0^{\cdot\wedge \tau_j(\mathcal{A})} \beta_j(r)d\vf_{0, r}(a_j) \right\rangle_s,
\end{align*}
thus
\begin{align*}
 \left| E\vf_{0, t}(0) \mathcal{E}^{|\mathcal{A}|}_{0, t}(\mathcal{A}; \beta) \right| & = \left|E\sum_{j=1}^{M} \int\limits_{0}^{t\wedge \tau_j(\mathcal{A})} \mathcal{E}^{|\mathcal{A}|}_{0, s}(\mathcal{A}; \beta) \beta_j(s)\1_{\{\vf_{0, s}(0) = \vf_{0, s}(a_j)\}}ds\right| 
\\
& \le 
\alpha^{\star}\sum_{j=1}^{M} \int\limits_{0}^{t} E \mathcal{E}^{|\mathcal{A}|}_{0, s}(\mathcal{A}; \beta) \1_\{\vf_{0, s}(0) = \vf_{0, s}(a_j)\}ds.
\end{align*} 

Suppose $a_j\not=0,$ the opposite case being of no interest. Since a process $\frac{\vf_{0, \cdot}(0) - \vf_{0, \cdot}(a_j)}{\sqrt{2}}$ is a Brownian motion stopped after hitting the origin the reflection principle~\cite{ito} gives:
$$
P\{\vf_{0, s}(0) = \vf_{0, s}(a_j)\} = \int_{0}^{\frac{2s}{a_j^2}}q(x)dx,
$$
where $q(x) = \frac{1}{\sqrt{2\pi x^3}}\e^{-\frac{1}{2x}}, x> 0.$ Put $q^{\star}=\max_{x>0} q(x)< +\infty.$ Obviously, $ \int_0^{\infty}q(x)dx=1.$ Take $\kappa\in(0;\frac{1}{2}).$ Using~(\ref{exp:mean}) and substituting $n^{-1}$ for $t:$
\begin{align*}
&  \left| E\vf_{0, t}(0) \mathcal{E}^{|\mathcal{A}|}_{0, t}(\mathcal{A}; \beta) \right| \le 
\alpha^{\star} \sum_{j=1}^{M} \int\limits_{0}^{t}  \left[ \left(E\mathcal{E}^{|\mathcal{A}|}_{0, s}(\mathcal{A}; \beta)\right)^2 \right]^{\frac{1}{2}} \left(\int\limits_0^{\frac{2s}{a_j^2}} q(x)dx\right)^{\frac{1}{2}} ds
\\ 
& \le
\alpha^{\star} \sum_{j=1}^{M} \int\limits_{0}^{t}  \left[ E\mathcal{E}^{|\mathcal{A}|}_{0, s}(\mathcal{A}; 2\cdot \beta)\exp\{t(\alpha^{\star})^2 M\} \right]^{\frac{1}{2}} \left(\int\limits_0^{\frac{2s}{a_j^2}} q(x)dx\right)^{\frac{1}{2}} ds
\\
& \le
\alpha^{\star} \e^{(\alpha^{\star})^2\frac{N}{2}} \sum_{j=1}^{M} \int\limits_{0}^{t}  \left(\int\limits_0^{\frac{2s}{a_j^2}} q(x)dx\right)^{\frac{1}{2}} \left( \1_{\{\frac{2s}{a_j^2}\le n^{-\kappa}\}} + \1_{\{\frac{2s}{a_j^2}> n^{-\kappa}\}} \right)ds
\\
& \le
\alpha^{\star} \e^{(\alpha^{\star})^2\frac{N}{2}} \sum_{j=1}^{M}  \int\limits_{0}^{t} \left( (q^{\star} n^{-\kappa})^{\frac{1}{2}}\1_{\{\frac{2s}{a_j^2}\le n^{-\kappa}\}} + \1_{\{\frac{2s}{a_j^2}> n^{-\kappa}\}}  \right) ds 
\\
& \le
\alpha^{\star} \e^{(\alpha^{\star})^2\frac{N}{2}} \left(   N(q^{\star})^{\frac{1}{2}} n^{-1-\frac{\kappa}{2}}  + \sum_{j=1}^{M} \int\limits_{0}^{t} \1_{\{ a_j^2 < 2sn^{\kappa}\}} ds \right).
\end{align*}
Returning to (\ref{app:main}) and ommiting the term for $k=0$ as an irrelevant one we have:
\begin{align*}
 \sum_{k=1}^{n-1} & \left| E \mathcal{E}^N_{0, \frac{k}{n}}(\mathcal{U}; \alpha) E\left(\xi_{n,k}\mathcal{E}^{|\mathcal{A}|}_{\frac{k}{n}, \frac{k+1}{n}}(\mathcal{A}; \alpha(\mathcal{K}))  \right) \Big|_{\begin{subarray}{l}{\mathcal{A}=\mathcal{U}_k,}\\
{\mathcal{K}=\mathcal{L}_k}
\end{subarray}
} \right| 
\\
& \le
\alpha^{\star} \e^{(\alpha^{\star})^2\frac{N}{2}} \sum_{k=1}^{n-1} E \mathcal{E}^N_{0, \frac{k}{n}}(\mathcal{U}; \alpha)   \left( N(q^{\star})^{\frac{1}{2}} n^{-1-\frac{\kappa}{2}}  +  \sum_{j=1}^{M}  \1_{\{ a_j^2 < 2sn^{\kappa}\}} ds \right) \Big|_{\begin{subarray}{l}{a_j = u_{k, j},}\\ {j = \overline{1, M_k}, }\\
{M= M_k}
\end{subarray}
}
\\
& \le
\alpha^{\star} \e^{(\alpha^{\star})^2\frac{N}{2}} \left( N (q^{\star})^{\frac{1}{2}} n^{-\frac{\kappa}{2}} + \sum_{k=1}^{n-1} E \mathcal{E}^N_{0, \frac{k}{n}}(\mathcal{U}; \alpha) \sum_{j=1}^{N} \int\limits_0^{n^{-1}} \1_{\{ \vf_{0, \frac{k}{n}}(u_j)^2 < 2sn^{\kappa}\}} ds \right) 
\\
& \le
\alpha^{\star} \e^{(\alpha^{\star})^2\frac{N}{2}} \left( N (q^{\star})^{\frac{1}{2}} n^{-\frac{\kappa}{2}} \right.
\\
& \phantom{abc} + \left. \sum_{k=1}^{n-1} \sum_{j=1}^{N}  \int\limits_0^{n^{-1}} \left[P\left\{ \vf_{0, \frac{k}{n}}(u_j)^2 < 2sn^{\kappa}\right\} E \left(\mathcal{E}^N_{0, \frac{k}{n}}(\mathcal{U}; \alpha)\right)^2  \right]^{\frac{1}{2}} ds \right)
\\
& \le
\alpha^{\star} \e^{(\alpha^{\star})^2\frac{N}{2}} \left( N (q^{\star})^{\frac{1}{2}} n^{-\frac{\kappa}{2}} \right.
\\
& \phantom{abc} + \left. \e^{(\alpha^{\star})^2\frac{N}{2}} \sum_{k=1}^{n-1} \sum_{j=1}^{N} \int\limits_0^{n^{-1}} \left[  P\{ \vf_{0, \frac{k}{n}}(u_j)^2 < 2sn^{\kappa}\} \right]^{\frac{1}{2}} ds  \right) 
\\
& \le
\alpha^{\star} \e^{(\alpha^{\star})^2\frac{N}{2}} \left( N (q^{\star})^{\frac{1}{2}} n^{-\frac{\kappa}{2}} + \e^{(\alpha^{\star})^2\frac{N}{2}} n^{-1} \sum_{k=1}^{n-1} \sum_{j=1}^{N}  P\{ \vf_{0, \frac{k}{n}}(u_j)^2 < 2n^{\kappa-1}\} \right)
\\
& \le
\alpha^{\star} \e^{(\alpha^{\star})^2\frac{N}{2}} \left( N (q^{\star})^{\frac{1}{2}} n^{-\frac{\kappa}{2}} \right. 
\\
& \phantom{abc} + \e^{(\alpha^{\star})^2\frac{N}{2}} n^{-1} \sum_{k=1}^{n-1} \sum_{j=1}^{N}  P\left\{ \left( \left(\frac{n}{k}\right)^{\frac{1}{2}}u_j + \mathcal{N}(0, 1)\right)^2 < 2n^{\kappa}k^{-1}\right\} 
\\
& \phantom{aaaaaaaaaaaaaaaa} \left. \cdot \left(\1_{\{k \le n^{2\kappa}\}} + \1_{\{k > n^{2\kappa}\}} \right) \right)
\\
& \le 
\alpha^{\star} \e^{(\alpha^{\star})^2\frac{N}{2}} \Bigg( N (q^{\star})^{\frac{1}{2}} n^{-\frac{\kappa}{2}} + \e^{(\alpha^{\star})^2\frac{N}{2}} N n^{-1} (n^{2\kappa} +1) 
\\
& \phantom{abc} +
\e^{(\alpha^{\star})^2\frac{N}{2}}  \sum_{j=1}^{N}  \max_{k=\overline{1, n-1}} P\left\{ \left( \left(\frac{n}{k}\right)^{\frac{1}{2}}u_j + 
\mathcal{N}(0, 1)\right)^2 < 2n^{-\kappa}\right\}\Bigg)
\\
& \le
\alpha^{\star} N\e^{(\alpha^{\star})^2\frac{N}{2}} \left((q^{\star})^{\frac{1}{2}} n^{-\frac{\kappa}{2}} + \e^{(\alpha^{\star})^2\frac{N}{2}} n^{-1} (n^{2\kappa} +1) \right.
\\
& \phantom{aaaaaaaaaaaaaaaaaaaaaa} \left. + \e^{(\alpha^{\star})^2\frac{N}{2}}   P\left\{ \mathcal{N}(0, 1)^2 < 2n^{-\kappa}\right\}\right).
\end{align*}
Thus we conclude.  
\qed

\end{document}